\documentclass[a4paper,12pt]{amsart}

\usepackage{amsmath,epsfig,amssymb,amsbsy,amscd,amsfonts,amstext,color,bm}

\usepackage[all]{xy}
\usepackage{bm}
\usepackage{comment}
\usepackage{enumerate}
\usepackage{circledsteps}
\usepackage{geometry}
\usepackage{lineno}

\makeatletter
\@namedef{subjclassname@2020}{%
  \textup{2020} Mathematics Subject Classification}
\makeatother

\usepackage[dvipdfmx,svgnames]{xcolor}
\usepackage{tikz}
\tikzset{v/.style={circle, draw, inner sep=2pt, minimum size=6pt, fill=white}}
\usetikzlibrary{patterns}

\theoremstyle{plain}
\newtheorem{theorem}{Theorem}[section]
\newtheorem*{theorem*}{Theorem}
\newtheorem*{theoremA*}{Theorem A}
\newtheorem*{theoremB*}{Theorem B}
\newtheorem{corollary}[theorem]{Corollary}

\theoremstyle{definition}
\newtheorem{definition}[theorem]{Definition}
\newtheorem{remark}[theorem]{Remark}

\newtheorem{proposition}[theorem]{Proposition}

\newcommand{\Ker}{\operatorname{Ker}}

\newcommand{\Tor}{\operatorname{Tor}}

\newcommand{\bC}{\mathbb{C}}

\newcommand{\bF}{\mathbb{F}}

\newcommand{\bQ}{\mathbb{Q}}

\newcommand{\bZ}{\mathbb{Z}}

\newcommand{\cA}{\mathcal{A}}

\def\qed{\hfill $\Box$}

\title{Alexander modules and Milnor fibers of topological line arrangements}
\author{Sakumi Sugawara}
\date{\today}
\address{Department of Mathematics, Faculty of Science, Hokkaido University, North 10, West 8, Kita-ku, Sapporo 060-0810, JAPAN. }
\email{sugawaras@math.sci.hokudai.ac.jp}

\begin{document}
\begin{abstract}
We study the Alexander module and a distinguished cyclic cover, which we call the \textit{topological Milnor fiber}, of the complements of topological line arrangements. For every topologically realizable combinatorial type with $n \geq 2$ lines, where $n$ is not a prime power, we prove that the first Betti numbers of its topological Milnor fibers are unbounded. Moreover, any finite abelian group occurs as a direct summand of the integral homology of a topological Milnor fiber within the given combinatorial type. As an application, we show that every topologically realizable combinatorial type with $n \geq 2$ lines admits a realization whose complement has no $2$-dimensional CW model.
\end{abstract}

\maketitle


\section{Introduction}
A finite set of lines in the complex projective plane is called a complex projective line arrangement. 
A central problem in the topology of arrangements is to determine which topological invariants are determined by the intersection poset \cite{orl-ter, dim-hyp}.
Whether the first Betti number of the Milnor fiber is combinatorially determined is a longstanding problem \cite{pap-suc-17}. 
Since the Milnor fiber is a cyclic cover of the complement, its homology can be studied using rank-one local systems and Alexander modules. 
Non-vanishing of a non-trivial monodromy eigenspace imposes a strong restriction on the combinatorics \cite{lib-eig,yos-mil}.
Recently, Xie-Yu-Zheng \cite{xyz} proved that, for any essential line arrangement, all monodromy eigenspaces corresponding to eigenvalues of order greater than $5$ always vanish. 

Torsion in the integral homology group of the Milnor fiber is another important aspect of this problem. It is an open problem whether the torsion in the homology group of the Milnor fiber is combinatorially determined. To the best of the author's knowledge, the icosidodecahedral arrangement introduced in \cite{yos-dou} is the only known reduced line arrangement whose Milnor fiber has torsion in its first integral homology (see also \cite{cds} for multiarrangements and \cite{den-suc-multi} for higher degree). 
The results below show that topological line arrangements admit much greater variation within a fixed combinatorial type.

A topological line arrangement $\cA$ is an arrangement of embedded $2$-spheres in the complex projective plane which topologically generalizes complex line arrangements (Section \ref{sec:pre}, see also \cite{rub-sta}). The \textit{topological Milnor fiber} $F_{\cA}$ is defined as a certain finite cyclic cover of the complement. If $\cA$ is a complex line arrangement, then $F_{\cA}$ coincides with the usual Milnor fiber.
The first result is on the first Betti numbers of the topological Milnor fibers.

\begin{theorem}\label{thm:main_mf}
Let $\cA$ be a topological line arrangement with $n \geq 2$ lines and set 
\[
D(n) = \{m \in \bZ_{>1} \mid m |n, \mbox{and $m$ is not a prime power}\}.
\]
Then, for any family $(k_m)_{m \in D(n)}$ of non-negative integers, there exists a topological line arrangement $\cA'$ such that $L(\cA) \cong L(\cA')$ and the topological Milnor fiber satisfies
\[
H_1 (F_{\cA'}; \bC) \cong H_1 (F_\cA; \bC) \oplus \bigoplus_{m \in D(n)} (\Lambda/ \Phi_{m}(t))^{\oplus k_m},
\]
as $\Lambda = \bC[t^{\pm 1}]$-modules, where $\Phi_m (t)$ is the $m$-th cyclotomic polynomial.
\end{theorem} 

\begin{corollary}
Let $e_{m} (\cA)$ be the dimension of the eigenspace corresponding to a primitive $m$-th root of unity. Then, for every $m \in D(n)$, 
\[
e_{m}(\cA') = e_{m}(\cA) + k_m.
\]
The first Betti number satisfies 
\[
b_1 (F_{\cA'}) = b_1 (F_\cA) + \sum_{m \in D(n)} \varphi (m)k_m,
\]
where $\varphi$ is Euler's totient function.
\end{corollary}
In particular, the first Betti number is unbounded if $n$ is not a prime power among realizations of the fixed combinatorial type $L(\cA)$.
The proof is based on the connected sum with $2$-knots, which is also used in \cite[Section 5]{sug-top} (see also \cite[p.~845]{lib-alex-82}, \cite[Section 2]{lib-alex-83}). 

The same local knotting construction also produces prescribed finite direct summands in the integral homology of the topological Milnor fiber.

\begin{theorem}\label{thm:mf_tor}
Let $\cA$ be a topological line arrangement with $n \geq 2$ lines. 
\begin{enumerate}[(1)]
\item Suppose $n$ is not a prime power. Then, for any finite abelian group $G$, there exists a topological line arrangement $\cA'$ such that $L(\cA) \cong L(\cA')$ and $H_i (F_{\cA'}; \bZ)$ contains $G$ as a direct summand for $i=1,2$.
\item Suppose $n=p^s$ is a prime power. Then, for any finite abelian group $G$ of order coprime to $ p$, there exists a topological line arrangement $\cA'$ such that $L(\cA) \cong L(\cA')$ and $H_i (F_{\cA'}; \bZ)$ contains $G$ as a direct summand for $i=1,2$.
\end{enumerate}
\end{theorem}

The restrictions on eigenvalue orders in Theorem \ref{thm:main_mf} and on the order of $G$ in Theorem \ref{thm:mf_tor} cannot be relaxed within the local knotting construction; see Remarks \ref{rmk:eigen} and \ref{rmk:tor}.

Moreover, the second Alexander module gives an obstruction to the existence of a $2$-dimensional CW model for the complement. 

\begin{theorem}\label{thm:2-dim}
Let $\cA$ be a topological line arrangement with $n \geq 2$ lines. Then, there exists a topological line arrangement $\cA'$ such that $L(\cA) \cong L(\cA')$ and $U(\cA')$ does not have the homotopy type of a $2$-dimensional CW complex.
\end{theorem}

Since $H_i (U)=0$ for $ i \geq 3$ for any topological line arrangement, this implies that $U(\cA')$ does not have the homotopy type of a minimal CW complex, yielding an alternative proof of \cite[Theorem 1.3]{sug-top}. 
By \cite[Theorem 1.2]{sug-top}, the complement of every symplectic line arrangement is minimal. Thus, if the $2$-knot used in the local knotting construction has a non-trivial first Alexander module over $\bQ$, the resulting arrangement cannot be symplectic (see Remark \ref{rmk:non-min}). Symplectic line arrangements share other topological properties with complex line arrangements, including divisibility relations for their Alexander polynomials \cite{awa-gol}. 
These parallels motivate a closer study of the Alexander modules and topological Milnor fibers of symplectic line arrangements. More broadly, it would be interesting to investigate the topology of topological line arrangements by methods that do not rely on connected sums with $2$-knots.

\vspace{3mm}
\textbf{Acknowledgements.}
The author would like to thank Anatoly Libgober for his inspiring comments and Masahiko Yoshinaga for his comments on torsion in homology.
This work is supported by JSPS KAKENHI Grant Number 25K23326.

\vspace{3mm}
\textbf{Declaration of AI use.}
The author used ChatGPT 6 Astra (OpenAI) during the research and preparation of this paper, particularly to discuss the proof of Theorem \ref{thm:mf_tor}, assist in checking mathematical arguments, search for relevant literature, and improve the English and exposition. The author has independently verified the mathematical arguments and references and takes full responsibility for the entire content of the paper.

\section{Preliminaries}\label{sec:pre}

\begin{definition}
A \textit{topological line arrangement} $\cA=\{H_0,\dots, H_{n-1}\}$ is a collection of locally flat embedded $2$-spheres in the complex projective plane such that each $H_{i}$ satisfies $[H_i] = [\bC P^1]\in H_2 (\bC P^2; \bZ)$, and two different lines intersect at a single point transversely and positively. 
Let $U=U(\cA)$ be the complement.
\end{definition}

\begin{definition}
The \textit{combinatorial type} (or the \textit{intersection poset}) is defined by the set of all intersections of lines, with the reverse inclusion ordering. We denote this poset by $L(\cA)$.
\end{definition}

From now on, we assume $n \geq 2$. In \cite{sug-top}, Sugawara proved that the first homology group is expressed as 
\[
H_1 (U; \bZ) \cong \bZ \langle m_0, \dots, m_{n-1} \rangle/ \langle m_0 +\dots + m_{n-1} \rangle \cong \bZ^{n-1},
\]
where $m_{i}$ denotes a positive meridian around $H_i$. Moreover, the cohomology ring is isomorphic to the projective Orlik-Solomon algebra of the combinatorial type.

The surjective homomorphism $\delta: \pi_1 (U) \to \bZ/ n \bZ, \delta (m_{i}) = 1$ defines an $n$-fold cyclic cover $F=F_{\cA} \to U$. We call this covering space the 
\textit{topological Milnor fiber} of $\cA$. If $\cA$ is a complex line arrangement, then this coincides with the usual Milnor fiber of the arrangement.

Also, by fixing the line $H_{0}$ (as the ``line at infinity''), we consider the homomorphism 
$\nu = \nu_{H_0}: \pi_1 (U) \to \bZ$, $\nu(m_i) = 1$ ($i\neq 0$), $\nu (m_0) = 1-n$. 
This is also surjective and is an integral lift of $\delta$. Let us consider the corresponding infinite cyclic cover $U^{\nu} \to U$. 
The deck transformation group $\bZ \cong \langle t \rangle $ acts on the $i$-th homology group $H_{i} (U^{\nu}; R)$; thus it is regarded as a module over $\Lambda_{R} = R[t^{\pm 1}]$, where $R$ is either $\bZ$ or a field. This $\Lambda_{R}$-module $H_{i}(U^{\nu}; R)$ is called the \textit{$i$-th Alexander module} of the topological line arrangement (associated with the distinguished line $H_0$).

\begin{proposition}\label{prop:homol_mf}
Let $M_{i} = \Ker (t^n-1: H_{i} (U^\nu; R) \to H_{i}(U^\nu; R)) $.
We have the following short exact sequence
\begin{eqnarray*}
0 \to H_i (U^{\nu}; R) / (t^n-1)H_i (U^{\nu}; R)  \to H_i (F; R) 
\to M_{i-1} \to 0,
\end{eqnarray*}
for $i\geq 1$ which splits when $i=1$. In particular, we have  an isomorphism of $\Lambda_R$-modules:
\[
H_1 (F; R) \cong H_1 (U^{\nu}; R) / (t^{n}-1) H_1 (U^{\nu}; R) \oplus R.
\]
\end{proposition}

\proof
Since $\delta$ is the reduction of $\nu$ modulo $n$, $F = U^{\nu}/\langle t^{n} \rangle$. We have the short exact sequence of cellular chain complexes
\[
0 \to C_{*} (U^\nu; R) \xrightarrow{t^n -1} C_* (U^\nu; R) 
\to C_* (F; R) \to 0.
\]
The long exact sequence in homology yields an exact sequence of $\Lambda_{R}$-modules
\[
0 \to H_{i}(U^\nu; R) / (t^{n}-1) H_{i}(U^\nu; R) \to H_i (F; R) 
\to M_{i-1} \to 0,
\]
for each $i \geq 1$. When $i=1$, since $H_0 (U^\nu; R) \cong R$ and $t$ acts on this trivially, we have 
\[
0 \to H_{1}(U^\nu; R) / (t^{n}-1) H_{1}(U^\nu; R) \to H_1 (F; R) 
\to R \to 0.
\]
Choose a loop $\mu \in \pi_1 (U)$ with $\nu (\mu) =1$. Then, the lift of $\mu^{n}$ defines a closed loop in $F$. Let $v \in H_1 (F; R)$ be its homology class. Note that $v$ is invariant under the deck transformation. 
By the definition of the connecting homomorphism, the map $R \to H_1 (F; R)$ sending the generator to $v$ defines a section.
Thus the above exact sequence splits, and the proposition follows.
\endproof

Now, we additionally suppose $R$ is algebraically closed and of characteristic zero and focus on the first homology group.
Since $t^n=1$, the deck transformation action on $H_1 (F; R)$ is diagonalizable, with eigenvalues among the $n$-th roots of unity. 
For each eigenvalue $\lambda \neq 1$, the eigenspace satisfies
\[
H_1 (F; R)_{\lambda} \cong H_1 (U^{\nu};R) \otimes \Lambda_R/ (t-\lambda)\Lambda_R. 
\]
Since the action is defined over $\bQ$, the dimension of the $\lambda$-eigenspace depends only on the order $m$ of $\lambda$ (see also \cite[Proposition 1.1]{coh-suc-mil}). We denote its dimension by $e_{m}(\cA)$. We also have 
\[
\dim_{R} H_1 (F; R)_{1} = b_1 (U) = n-1.
\]
Therefore, the first Betti number of the topological Milnor fiber is expressed as
\[
b_1 (F) = n-1 + \sum_{1<m|n} \varphi(m)e_{m}(\cA).
\]

\section{Proof of the results}
In this section, we prove the main results. We first recall a realization theorem for Alexander polynomials of $2$-knots. For a $2$-knot $K$ with the complement $X_K =S^4 \setminus K$, since $H_1 (X_K; \bZ) \cong \bZ$, there exists a unique homomorphism $\rho: \pi_1 (X_K) \to \bZ$ sending the positive meridian to the generator.
The $i$-th homology group $H_i (X^{\rho}_K; R)$ of the corresponding infinite cyclic cover is called the $i$-th Alexander module of $K$, and the order $\Delta_K (t)$ of $H_1 (X_K^{\rho}; \bQ)$ as a $\Lambda_{\bQ}$-module is called the Alexander polynomial of $K$.
For the Alexander polynomial of a $2$-knot, the following realization theorem is known.

\begin{proposition}[\cite{kin}]\label{prop:kin}
For any integer polynomial $P(t)$ with $P(1) = 1$, there exists a $2$-knot $K$ such that $\Delta_K (t) = P(t)$.
\end{proposition}
In particular, if $P(t)$ is square-free, the structure theorem implies that  $H_1 (X^{\rho}_K; R) \cong \Lambda_{R}/ (P(t))$ for every field $R$ of characteristic zero.

From a $2$-knot $K$ and a topological line arrangement $\cA = \{H_0,\dots, H_{n-1}\}$ with $H_0$ fixed as the line at infinity, we can obtain a new topological line arrangement $\cA'$ with $L(\cA') \cong L(\cA)$, by taking a connected sum of a line $H (\neq H_0)$ and the $2$-knot $K$ at a smooth point of $H$.
Let $U'$ be the complement of the resulting arrangement and $\nu': \pi_1 (U') \to \bZ$ be the induced homomorphism.
In the connected sum decomposition, the gluing region is homotopy equivalent to a meridian circle. The character $\nu'$ has value $1$ on this circle, so its lift to the infinite cyclic cover is connected and acyclic. The same is true of the standard disk-complement piece.
The Mayer-Vietoris argument yields the $\Lambda_R$-module isomorphism
\[
H_i((U')^{\nu'}; R) \cong H_i (X_K^{\rho}; R) \oplus H_i (U^{\nu}; R)
\]
for $i \geq 1$. Thus, by Proposition \ref{prop:homol_mf}, we have
\begin{equation}\label{eq:first}
H_1 (F_{\cA'}; R) \cong H_1 (F_{\cA}; R)  \oplus H_1 (X_K^{\rho}; R) / (t^n-1) H_1 (X_K^{\rho}; R).
\end{equation}
Similarly, by the Mayer-Vietoris argument for the $n$-fold cyclic cover, we obtain
\begin{equation}\label{eq:second}
H_2 (F_{\cA'}; R) \cong H_2 (F_{\cA}; R) \oplus H_2 (X_{K,n}; R)
\end{equation}
where $X_{K,n}$ is the $n$-fold cyclic cover of $X_{K}$. 

Now, we are ready to prove the main results.

\vspace{3mm}

\textit{Proof of Theorem \ref{thm:main_mf}.} 
In this proof, we consider the Alexander module over $\bC$.
For $m \in D(n)$, the cyclotomic polynomial satisfies $\Phi_m (1)=1$. 
Thus, by the observation immediately following Proposition \ref{prop:kin}, there exists a $2$-knot $K_m$ whose Alexander module is isomorphic to $\Lambda_\bC/ (\Phi_m (t))$. Then, the Alexander module of the connected sum $K = \sharp_{m \in D(n)}(\sharp_{k_m} K_m)$ is 
$\bigoplus_{m \in D(n)}(\Lambda/ \Phi_m (t))^{\oplus k_m}$. 
Let $\cA'$ be the topological line arrangement obtained by taking the connected sum of a line other than $H_0$ with $K$. 
Since $\Phi_m (t)$ divides $t^n-1$ for every $m \in D(n)$, multiplication by $t^n-1$ vanishes on $H_1 (X_K^{\rho};\bC)$. Therefore, by (\ref{eq:first}),
\begin{eqnarray*}
H_1(F_{\mathcal A'};\mathbb C)
&\cong& H_1(F_{\mathcal A};\mathbb C)
  \oplus H_1 (X_K^\rho; \bC)/(t^n-1)H_1 (X_K^\rho;\bC) \\
&\cong& H_1(F_{\mathcal A};\mathbb C)
  \oplus \bigoplus_{m\in D(n)}
  \left ( \Lambda/\Phi_m(t) \right )^{\oplus k_m}
\end{eqnarray*}
as $\Lambda_{\bC}$-modules. This proves the theorem.
\qed

\begin{remark}\label{rmk:eigen}
For a prime power $m=p^r$, the identity $\Phi_m(1) = p$ prevents $\Phi_m (t)$ from dividing an integral knot Alexander polynomial with value $\pm 1$ at 1. Thus, local knotting adds no eigenspaces of prime-power order.
More generally, the Papadima-Suciu inequality \cite[Theorem C]{pap-suc-10} states that the dimension of the corresponding monodromy eigenspace is bounded by the mod-$p$ Aomoto-Betti number, which is combinatorially determined.
\end{remark}

Next, we consider the torsion in the homology group of the topological Milnor fiber. 
The local knotting realizes prescribed finite abelian groups as direct summands.

\vspace{3mm}

\textit{Proof of Theorem \ref{thm:mf_tor}}.
By the local knotting formula above, it suffices to consider the homology group of an unbranched $n$-fold cyclic cover of the complement of a  $2$-knot.
Let $q$ be a prime factor of $n$ and $k \geq 2$ be an integer coprime to $q$. Let $K_{k,q}$ be the $2$-knot obtained from the torus knot $T(k,q)$ by $q$-twist spinning. 
By \cite{zee}, the complement $X_{K_{k,q}}$ fibers over $S^1$ with the fiber $\Sigma_{q}(T(k,q)) \setminus \{x\}$, where $\Sigma_{q}(T(k,q))$ is the $q$-fold branched cyclic cover of $S^3$ along $T(k,q)$ and $x$ is a point on the ramification locus. 
Its monodromy map of the bundle is known to be the deck transformation of the branched cover. 
Thus, the Alexander module of $K_{k,q}$ satisfies 
\[
H_i (X_{K_{k,q}}^{\rho}; \bZ) \cong H_i (\Sigma_{q}(T(k,q)); \bZ)
\]
for $i =1,2$ and the multiplication by $t$ is induced by a deck transformation of the branched cover. In particular, $t^{q}$ acts trivially on $H_i (X_{K_{k,q}}^{\rho}; \bZ)$. Since $q$ is a factor of $n$, by an argument similar to the proof of Proposition \ref{prop:homol_mf},
the homology group of the $n$-fold cyclic cover $X_{K_{k,q}, n}$ satisfies 
\[
H_1 (X_{K_{k,q}, n} ;\bZ) \cong  \bZ \oplus H_1 (X_{K_{k,q}}^{\rho}; \bZ) / (t^{n}-1) H_1 (X_{K_{k,q}}^{\rho}; \bZ) \cong \bZ \oplus H_1 (X_{K_{k,q}}^{\rho}; \bZ). 
\]
By a direct computation (see also \cite[Corollary D]{nof-wil} with $(r,s,n)=(k,q,q)$), we have that 
$H_1 (X_{K_{k,q}}^{\rho}; \bZ) \cong (\bZ/ k\bZ)^{q-1}$. 
Also, since $H_2 (X_{K_{k,q}}^{\rho}; \bZ) =0$ by Poincar\'e duality, we have
\[
H_2 (X_{K_{k,q}, n};\bZ) \cong \Ker(t^{n}-1: H_1 (X_{K_{k,q}}^{\rho} ;\bZ) \to H_1 (X_{K_{k,q}}^{\rho} ;\bZ)) \cong H_1 (X_{K_{k,q}}^{\rho} ;\bZ)
\]
by the first exact sequence of Proposition \ref{prop:homol_mf}. Thus local knotting with $K_{k,q}$ adds a direct summand $(\bZ/ k\bZ)^{q-1}$ in each of degrees $1$ and $2$ by (\ref{eq:first}) and (\ref{eq:second}).

Let $G$ be any finite abelian group. We express $G$ as $G \cong \bigoplus_{i=1}^{r} (\bZ/\ell_{i}^{a_i} \bZ)$ where $\ell_i$ is prime and $a_i \geq 1$.
Suppose $n$ is not a prime power. Then, $n$ has at least two distinct prime factors, so for every $i$, we can choose a prime $q_{i} | n$ with $q_i \neq \ell_i$. Set $k_i = \ell_i^{a_i}$ and take $K=K_{k_1, q_1} \sharp \cdots \sharp K_{k_r, q_r}$. Then, the connected sum with $K$ adds $\bigoplus_{i=1}^{r} (\bZ/ k_{i} \bZ)^{q_i -1 }$ which contains $G$ as a direct summand. This proves (1).

If $n=p^s$ is a prime power and $|G|$ is coprime to $p$, then every $\ell_i \neq p$. Taking $q_i = p$ for all $i$ gives the same construction and proves (2).
\qed

\begin{remark}\label{rmk:tor}
Suppose $n = p^{s}$ is a prime power. Then, this knot connected sum adds no $p$-primary torsion for any $2$-knot. To see this, let $K$ be any $2$-knot and set $B_i = H_i (X_K^{\rho}; \bF_{p})$. 
We have the following short exact sequence of chain complexes
\[
0 \to C_{*}(X_{K}^{\rho}; \bF_p) \xrightarrow{t-1} C_{*} (X_{K}^{\rho}; \bF_p)
\to C_{*}(X_K ; \bF_p) \to 0.
\]
The long exact sequence shows that $(t-1)$ acts invertibly on each $B_i$ $(i \geq 1)$, and the same holds for $t^{n}-1 = (t-1)^n$. Hence $H_j (X_{K,n}; \bF_p) =0$ ($j \geq 2$). Therefore, by the universal coefficient theorem, 
$\Tor_1^{\bZ} (H_{j-1}(X_{K,n};\bZ), \bF_{p})$ vanishes for $j-1\geq1$.
Thus, both $H_1 (X_{K,n};\bZ)$ and $H_2 (X_{K,n}; \bZ)$ have no non-trivial $p$-torsion. 
\end{remark}

Finally, we prove Theorem \ref{thm:2-dim} as an application of the second Alexander module.

\vspace{3mm}

\textit{Proof of Theorem \ref{thm:2-dim}}.
Let $K$ be a spun trefoil. Then, its Alexander polynomial is $\Delta_{K}(t)=t^2 - t+1$ and its first Alexander module is $\Lambda_{\bQ} /(t^2-t+1)$. By \cite[Theorem 3.4]{lev}, the second Alexander module is isomorphic to $\Lambda_{\bQ} / (\Delta_K(t^{-1}))$. 
In particular, this is a torsion $\Lambda_{\bQ}$-module.
Therefore, the second Alexander module of the arrangement $\cA'$ obtained by local knotting with $K$ contains a non-trivial $\Lambda_{\bQ}$-torsion summand.

Now we suppose that $U'$ has the homotopy type of a $2$-dimensional CW complex; then so does $(U')^{\nu'}$. Thus, 
\[
H_2 ((U')^{\nu'}; \bQ) = \Ker (\partial_2: C_2 ((U')^{\nu'}; \bQ) \to C_1 ((U')^{\nu'}; \bQ))
\]
is a free $\Lambda_{\bQ}$-module. This is a contradiction. Thus $U'$ is not homotopy equivalent to a $2$-dimensional CW complex.
\qed

\begin{remark}\label{rmk:non-min}
Theorem \ref{thm:2-dim} also follows directly from
Theorem \ref{thm:mf_tor}, since every covering space of a $2$-dimensional CW complex has torsion-free second integral homology. 
By \cite[Theorem 3.4 (ii)]{lev}, if a $2$-knot $K$ satisfies $H_1 (X_K^{\rho}; \bQ) \neq 0$, then $H_2 (X_K^{\rho}; \bQ)$ is a nonzero torsion $\Lambda_\bQ$-module.
Thus, taking the connected sum with such a knot yields an arrangement whose complement 
does not have the homotopy type of a $2$-dimensional CW complex.
\end{remark}


\begin{thebibliography}{999999}


\bibitem[AG24]{awa-gol}
H.~Awada, M.~Golla,
\emph{Oka and Alexander polynomials of symplectic curves and divisibility relations},
arXiv:2412.15792.

\bibitem[CDS03]{cds}
D. C. Cohen, G. Denham, A. Suciu, 
\emph{Torsion in Milnor fiber homology},
Algebr. Geom. Topol. \textbf{3} (2003), 511--535. 

\bibitem[CS95]{coh-suc-mil}
D. C. Cohen, A. Suciu, 
\emph{On Milnor fibrations of arrangements},
J. Lond. Math. Soc. {\bf 51} (1995), no. 1, 105--119.

\bibitem[DS14]{den-suc-multi}
G. Denham, A. Suciu, 
\emph{Multinets, parallel connections, and Milnor fibrations of arrangements}, Proc. Lond. Math. Soc. {\bf 108} (2014), no. 6, 1435--1470.

\bibitem[Di17]{dim-hyp}
A. Dimca, 
\emph{Hyperplane arrangements. An introduction},
Universitext, Springer, Cham, 2017, xii+200 pp.

\bibitem[Ki61]{kin}
S.~Kinoshita,
\emph{On the Alexander polynomials of 2-spheres in a 4-sphere},
Ann. of Math. 74 (1961), 518--531. 

\bibitem[Le77]{lev}
J. Levine,
\emph{Knot modules.~I},
Trans. Amer. Math. Soc., {\bf 229} (1977), 1--50.

\bibitem[Li82]{lib-alex-82}
A.~Libgober, 
\emph{Alexander polynomial of plane algebraic curves and cyclic multiple planes},
Duke Math. J., 49 (1982), 833--851.

\bibitem[Li83]{lib-alex-83}
A.~Libgober, 
\emph{Alexander modules of plane algebraic curves},
Low-dimensional topology (San Francisco, Calif., 1981), 231--247.
Contemp. Math., 20, American Mathematical Society, Providence, RI, 1983. 


\bibitem[Li02]{lib-eig}
A. Libgober, 
\emph{Eigenvalues for the monodromy of the Milnor fibers of arrangements}, 
Trends in singularities, 141--150, Trends Math., Birkh\"{a}user, Basel, 2002.

\bibitem[NW21]{nof-wil}
V. Noferini, G. Williams,
\emph{Matrices in companion rings, Smith forms, and the homology of 3-dimensional Brieskorn manifolds},
 J. Algebra, {\bf 587} (2021), 1--19.
 
\bibitem[OT92]{orl-ter}
P. Orlik, H. Terao, 
\emph{Arrangements of Hyperplanes}, 
Grundlehren Math. Wiss. {\bf 300}, Springer-Verlag, Berlin, 1992. 

\bibitem[PS10]{pap-suc-10}
S. Papadima, A. Suciu,
\emph{The spectral sequence of an equivariant chain complex and homology with
local coefficients}, 
Trans. Amer. Math. Soc. {\bf 362} (2010), no. 5,
2685--2721.

\bibitem[PS17]{pap-suc-17}
S.~Papadima, A.~Suciu,
\emph{The Milnor fibration of a hyperplane arrangement: from modular resonance to algebraic monodromy},
Proc. Lond. Math. Soc., 114 (2017), 961--1004.

\bibitem[RS19]{rub-sta}
D. Ruberman, L. Starkston,
\emph{Topological realizations of line arrangements},
Int. Math. Res. Not. IMRN (2019),
no. 8, 2295--2331.

\bibitem[Su26]{sug-top}
S. Sugawara, 
\emph{Topological line arrangements and their topological invariants}, arXiv:2607.23570.

\bibitem[XYZ26]{xyz}
B.~Xie, C. Yu, Z.~Zheng,
\emph{Monodromy eigenvalues of Milnor fibers for line arrangements},
arXiv:2609.06416.

\bibitem[Yo13]{yos-mil}  
M. Yoshinaga, 
\emph{Milnor fibers of real line arrangements}, 
Journal of Singularities, \textbf{7} (2013), 220--237.

\bibitem[Yo20]{yos-dou}
M. Yoshinaga, 
\emph{Double coverings of arrangement complements and 2-torsion in Milnor fiber homology}, 
Eur. J. Math. {\bf 6} (2020), no. 3, 1097--1109.

\bibitem[Ze65]{zee}
E. C. Zeeman,
\emph{Twisting spun knots},
Trans. Amer. Math. Soc. {\bf 115} (1965), 471--495.

\end{thebibliography}
\end{document}